\documentclass[12pt]{amsart}

\usepackage{amscd,amsmath,amssymb,amsthm,mathrsfs}
\usepackage{geometry}
\usepackage{shuffle}
\usepackage{mathtools}
\usepackage{hyperref}
\mathtoolsset{showonlyrefs}
\usepackage{tikz}
\usetikzlibrary{matrix,arrows,decorations.pathmorphing}
\usepackage{tikz-cd}

\numberwithin{equation}{section}

\def\bC{\mathbb C}
\def\bG{\mathbb G}
\def\bL{\mathbb L}
\def\bP{\mathbb P}
\def\bQ{\mathbb Q}

\def\bZ{\mathbb Z}

\def\cA{\mathcal A}

\def\cL{\mathcal L}

\def\cO{\mathcal O}

\def\cU{\mathcal U}

\def\sx{\mathsf{x}}

\DeclareMathOperator{\ad}{ad}

\DeclareMathOperator{\Hom}{Hom}

\DeclareMathOperator{\Int}{Int}

\DeclareMathOperator{\Pic}{Pic}

\DeclareMathOperator{\spn}{Span}

\theoremstyle{plain}
\newtheorem{prop}{Proposition}[section]
\newtheorem{thm}[prop]{Theorem}

\theoremstyle{definition}

\newtheorem{rmk}[prop]{Remark}
\newtheorem{exmp}[prop]{Example}

\title[Towards 
Algebraic iterated integrals for elliptic curves]{Towards algebraic iterated integrals for elliptic curves via the universal vectorial extension}
\author{Tiago J. Fonseca}
\address{Mathematical Institute, University of Oxford,
	Andrew Wiles Building, Radcliffe Observatory Quarter, Woodstock Road, Oxford OX2 6GG, United Kingdom}
\email{tiago.jardimdafonseca@maths.ox.ac.uk}

\author{Nils Matthes}
\address{Mathematical Institute, University of Oxford,
	Andrew Wiles Building, Radcliffe Observatory Quarter, Woodstock Road, Oxford OX2 6GG, United Kingdom}
\email{nils.matthes@maths.ox.ac.uk}

\subjclass[2010]{11F67 (11M32)}
\keywords{Elliptic curves, periods, iterated integrals, universal vectorial extension}

\begin{document}

\begin{abstract}
For an elliptic curve $E$ defined over a field $k \subset \bC$, we study iterated path integrals of logarithmic differential forms on $E^{\dagger}$, the universal vectorial extension of $E$. These are generalizations of the classical periods and quasi-periods of $E$, and are closely related to multiple elliptic polylogarithms and elliptic multiple zeta values. Moreover, if $k$ is a finite extension of $\bQ$, then these iterated integrals along paths between $k$-rational points are periods in the sense of Kontsevich--Zagier.
\end{abstract}
\maketitle

\section{Introduction}

This is a report on work in progress; details and further results will appear elsewhere.
The purpose of this note is merely to indicate an algebraic approach to the study of iterated path integrals on a once-punctured elliptic curve $\widetilde{E}=E\setminus \{O\}$. Our main motivation is to clarify the arithmetic structure of multiple elliptic polylogarithms, \cite{BrownLevin}, and their special values at torsion points, in particular putting them into the framework of (motivic) periods, \cite{Brown:NotesMotivicPeriods,KZ}.

Several approaches to iterated path integrals on $\widetilde{E}$ have already been studied in the literature, notably in \cite{BannaiKobayashiTsuji,BrownLevin,Enriquez:Emzv,LevinRacinet,MaLuo}. These were at least in part motivated by the study of the motivic fundamental group of $\bP^1 \setminus \{0,1,\infty\}$, \cite{Deligne:P1}, and its relation to mixed Tate motives over $\bZ$, \cite{DeligneGoncharov}. However, a major difficulty in the elliptic case is that the complex of global algebraic forms on $E$ with logarithmic poles along $O$ does not compute de Rham cohomology of $\widetilde{E}$, unlike the situation for $\bP^1 \setminus \{0,1,\infty\}$. To circumvent this, one can either allow higher order poles at $O$, \cite{BannaiKobayashiTsuji,LevinRacinet}, or non-algebraic forms, \cite{BrownLevin}, or work locally with a \v{C}ech covering of $E$, \cite{MaLuo}.
In this paper, we propose a new global algebraic approach via the universal vectorial extension of $E$.

\subsection{Algebraic invariants of manifolds}

We begin by giving some context for iterated integrals in general. Let $M$ be a connected complex manifold\footnote{For simplicity, we work with complex manifolds since this is the case we are ultimately interested in, although many results that we state are valid for real manifolds as well.} and denote by $\cA^{\bullet}(M)$ the $C^{\infty}$-de Rham complex, consisting of smooth differential forms on $M$ with values in $\bC$. The de Rham cohomology groups $H^n_{\rm dR}(M)$ are then by definition the cohomology groups $H^n(\cA^{\bullet}(M))$, and the classical de Rham theorem says that the integration map
\begin{equation}
\begin{aligned}
H^n_{\rm dR}(M) &\rightarrow H^n_{\rm sing}(M;\bZ)\otimes_\bZ \bC\\
[\omega] &\mapsto \left([\sigma] \mapsto \int_\sigma\omega\right),
\end{aligned}
\end{equation}
is an isomorphism of $\bC$-vector spaces, where $H^n_{\rm sing}(M;\bZ)$ denotes the singular cohomology of $M$. This result expresses a relation between certain analytic and topological data associated to $M$: singular cohomology classes with coefficients in $\bC$ can be computed using smooth differential forms.

The singular (co-) homology groups are only one kind of algebraic invariant that can be associated to $M$. Other examples include the homotopy groups $\pi_n(M,p)$, for $n\geq 1$, where $p \in M$ is a base point. Already the case $n=1$ is very interesting: $\pi_1(M,p)$ is the fundamental group of $M$ (based at $p$), and one has the Hurewicz map $\pi_1(M,p) \rightarrow H_1(M;\bZ)$, which induces an isomorphism $\pi_1(M,p)^{\rm ab}\cong H_1(M;\bZ)$. Since $\pi_1(M,p)$ is in general non-abelian, it is a finer invariant than $H_1(M;\bZ)$.

\subsection{Iterated path integrals and Chen's \texorpdfstring{$\pi_1$}{pi}-de Rham theorem}

In the 1970s, K.-T. Chen studied de Rham type theorems for homotopy groups. In particular, his $\pi_1$-de Rham theorem generalizes de Rham's theorem for $n=1$, \cite{Chen:IteratedPathIntegrals}; see also \cite{Hain:Bowdoin} for another exposition. We will assume from now on that $H_1(M;\bQ)$ is a finite-dimensional $\bQ$-vector space. For a point $p \in M$, let $P(M,p)$ be the set of all piecewise smooth loops $\gamma: [0,1] \rightarrow M$ based at $p$, and denote by $H^0(B_{\ell}(M,p)) \subset \Hom(P(M,p),\bC)$ the set of all linear combinations of iterated path integrals of length $\leq \ell$ which are homotopy functionals (see Section \ref{ssec:notationandconventions} for definitions). There is a well-defined $\bC$-linear map
\begin{equation}
\begin{aligned}
\Int^\ell_p: H^0(B_{\ell}(M,p)) &\rightarrow \Hom(\bZ\pi_1(M,p)/I^{\ell+1},\bC)\\
\sum_{n=0}^{\ell}\sum_{i_1,\ldots,i_n}\int\omega_{i_1}\ldots\omega_{i_n} &\mapsto \left( \gamma\mapsto \sum_{n=0}^{\ell}\sum_{i_1,\ldots,i_n}\int_\gamma\omega_{i_1}\ldots\omega_{i_n}  \right),
\end{aligned}
\end{equation}
where $\bZ\pi_1(M,p)$ is the group ring of the fundamental group, and $I:=\ker(\bZ\pi_1(M,p) \stackrel{\gamma\mapsto 1}\longrightarrow \bZ)$ is the augmentation ideal.
\begin{thm}[Chen] \label{thm:Chenintro}
The map $\Int^{\ell}_p$ is an isomorphism.
\end{thm}
Taking the limit over all $\ell \geq 0$, Theorem \ref{thm:Chenintro} induces an isomorphism of commutative filtered Hopf algebras $H^0(B(M,p)) \cong \varinjlim_\ell \Hom(\bZ\pi_1(M,p)/I^{\ell+1},\bC)$ where $H^0(B(M,p)):=\bigcup_{\ell\geq 0}H^0(B_\ell(M,p))$. There is also a variant for two different base points $p,q \in M$. The case $\ell=1$ recovers de Rham's theorem for $n=1$. Indeed, $\bZ\pi_1(M,p)/I^2 \cong \bZ\oplus I/I^2 \cong \bZ\oplus H_1(M;\bZ)$, by the Hurewicz isomorphism, and $H^0(B_1(M,p)) \cong \bC\oplus H^1_{\rm dR}(M)$, since a path integral $\int\omega$ is homotopy invariant if and only if $\omega \in \cA^1(M)$ is closed.

\subsection{The reduced bar complex}

Finding all homotopy invariant iterated path integrals is a rather subtle problem, for example there are double integrals $\int\omega_1\omega_2$, with $\omega_1,\omega_2$ both closed, which aren't homotopy functionals. Chen solved this problem using the reduced bar de Rham complex $B^{\bullet}(\cA^{\bullet}(M),p)$, \cite{Chen:ReducedBarComplex}, which is a differential graded $\bC$-algebra constructed from the smooth de Rham complex $\cA^{\bullet}(M)$ together with its augmentation $\varepsilon_p: \cA^{\bullet}(M) \rightarrow \bC$ given by evaluation at $p$.
\begin{thm}[Chen] \label{thm:ChenintroII}
There is a natural isomorphism of commutative filtered Hopf algebras: $H^0(B^{\bullet}(\cA^{\bullet}(M),p)) \stackrel{\sim}\rightarrow H^0(B(M,p))$.
\end{thm}
This theorem gives a purely algebraic description of all homotopy invariant iterated integrals on $M$. More generally, if $A^{\bullet} \subset \cA^{\bullet}(M)$ is a differential graded $\bC$-subalgebra which is quasi-isomorphic to $\cA^{\bullet}(M)$ and such that the restriction of $\varepsilon_p$ is non-trivial, we obtain an isomorphism $B^\bullet(A^{\bullet},p) \cong B^{\bullet}(\cA^\bullet(M),p)$ of the associated reduced bar complexes. In the special case where $A^{\bullet}$ is connected, i.e. $A^0 \cong \bC$, the reduced bar complex does not depend on the point $p$ and we write $B^{\bullet}(A^{\bullet})$ instead of $B^{\bullet}(A^{\bullet},p)$.

\subsection{Algebraic iterated path integrals and their periods}

Now assume that $M=X(\bC)$ where $X$ is a smooth algebraic variety which is the complement $X=\overline{X}\setminus D$ in a smooth projective variety $\overline{X}$ of a normal crossings divisor $D$, all defined over a subfield $k \subset \bC$.\footnote{Such an $M$ necessarily satisfies $\dim_\bQ H_1(M;\bQ)<\infty$.} In this situation, one has the differential graded $k$-subalgebra $H^0(\overline{X},\Omega^{\bullet}_{\overline X}\langle D\rangle) \subset \cA^{\bullet}(M)$ of global algebraic differential forms on $\overline{X}$ with logarithmic poles along $D$. In general, the canonical map 
\begin{equation} \label{eqn:embedding}
H^0(\overline{X},\Omega^{\bullet}_{\overline X}\langle D\rangle) \otimes_k\bC \hookrightarrow \cA^{\bullet}(M)
\end{equation}
is not a quasi-isomorphism, not even if $X$ is affine.\footnote{Instead, one needs to work with the full complex of sheaves 
$
\Omega^{\bullet}_{\overline{X}}\langle D\rangle$ on $\overline{X}$, whose \textit{hypercohomology groups} are functorially isomorphic to $H^\ast_{\rm dR}(M)$, \cite[Proposition 3.1.8]{Deligne:HodgeII} or \cite[Proposition 3.1.16]{HuberMuellerStach}.} For example, if $X$ is a curve, \eqref{eqn:embedding} is a quasi-isomorphism if and only if $H^1(\overline{X},\cO_{\overline{X}})=\{0\}$, i.e. $\overline{X}$ has genus zero.

However, already the case $\overline{X}=\bP^1_k$, and $D=\{p_1,\ldots,p_n,\infty\} \subset \overline{X}(k)$ is very interesting. In this case, we have
\begin{equation}
H^0(B(M,p))\cong \spn_k\left\{\int\omega_{i_1}\ldots\omega_{i_n} \, \Big\vert \, i_j \in \{p_1,\ldots,p_n\}\right\} \otimes_k \bC,
\end{equation}
for any base point $p \in M$, where $\omega_{p_i}=dz/(z-p_i) \in H^0(\overline{X},\Omega^1_{\overline{X}}\langle D\rangle)$, and $z$ is the standard coordinate on $\bP^1_k\setminus \{\infty\}$. This corresponds to a $k$-structure on $\varinjlim_\ell \Hom(\bZ\pi_1(X(\bC),p),\bC)$ under the integration map $\Int_p$, which is given by the values of the iterated integrals $\int_\gamma\omega_{i_1}\ldots\omega_{i_n} \in \bC$, for loops $\gamma$ based at $p$. If $k$ is a finite extension of $\bQ$ and $p\in X(k)$, then these complex numbers are periods in the sense of Kontsevich and Zagier, \cite{KZ}. Furthermore, since the forms $\omega_{p_i}$ have only logarithmic poles, the formalism of tangential base points applies and we may replace the base point $p$ by a non-zero tangent vector $\vec{v}_p$ at $p \in D$, \cite[\S 15]{Deligne:P1}. In the special case $k=\bQ$, $D=\{0,1,\infty\}$ and $\vec{v}_p=\pm\vec{1}_p$, we obtain $\bQ[2\pi i]$-linear combinations of multiple zeta values in this way.

\subsection{The case of elliptic curves}

Now consider an elliptic curve $E$ defined over a subfield $k\subset \bC$, with origin $O \in E(k)$, and let $\widetilde{E}:=E\setminus \{O\}$ be the once-punctured elliptic curve. In that case, the map \eqref{eqn:embedding} is not a quasi-isomorphism. However, following an idea already used in \cite{EnriquezEtingof} (attributed to Deligne), we get around this problem by replacing $E$ with its universal vectorial extension $E^{\dagger}$, \cite{MazurMessing,Rosenlicht}. This is a commutative $k$-group scheme which fits in a short exact sequence
\begin{equation}
0\longrightarrow \underline{\Omega}^1_E \longrightarrow E^{\dagger}\stackrel{\pi}\longrightarrow E\longrightarrow 0,
\end{equation}
where $\underline{\Omega}^1_E \cong \bG_a$ denotes the additive group scheme of global algebraic differentials on $E$. Let $\widetilde{E}^{\dagger}:=E^{\dagger}\setminus D$, where $D=\pi^{-1}(O)$. The key point is that, while the complex manifolds $\widetilde{E}^{\dagger}(\bC)$ and $\widetilde{E}(\bC)$ have quasi-isomorphic $C^{\infty}$-de Rham complexes, the natural map
\begin{equation}
H^0(E^{\dagger},\Omega^{\bullet}_{E^{\dagger}}\langle D\rangle) \otimes_k\bC \hookrightarrow \cA^{\bullet}(\widetilde{E}^{\dagger}(\bC)), 
\end{equation}
is a quasi-isomorphism, see Proposition \ref{prop:quasiisomorphism}. In particular, for any base point $p \in \widetilde{E}^{\dagger}(k)$, we get an explicit $k$-structure $H^0(B^{\bullet}(A^{\bullet}_{E^\dagger}))$ on $H^0(B(\widetilde{E}^{\dagger}(\bC),p))$ which induces, via Chen's $\pi_1$-de Rham theorem, a $k$-structure on $\varinjlim_\ell \Hom(\bZ\pi_1(\widetilde{E}^{\dagger}(\bC),p),\bC)$. Again, in the case where $k$ is finite over $\bQ$, this $k$-structure is given by periods in the sense of Kontsevich--Zagier, the simplest examples of which are
\begin{equation} \label{eqn:quasiperiods}
\int_{\gamma_i}\omega, \quad \int_{\gamma_i}\eta,
\end{equation}
where $\gamma_1,\gamma_2$ is a symplectic basis of $H_1(E(\bC);\bZ)$, and the classes of $\omega$, $\eta$ form a basis of $H^1_{\rm dR}(E)$, the algebraic de Rham cohomology of $E$. Classically, the numbers \eqref{eqn:quasiperiods} are known as the (quasi-) periods of $E$, \cite[\S 14.4]{HuberMuellerStach}.

One can also define elliptic analogues of multiple polylogarithms and multiple zeta values in this setting, which differ slightly from the ones introduced in \cite{BrownLevin,Enriquez:Emzv}. For this, one needs to replace the base point $p \in E^{\dagger}(k)$ with a tangential base point, and regularize the corresponding iterated path integrals, which is possible again since we are dealing with logarithmic poles rather than arbitrary ones. Furthermore, one also needs to extend such a formalism to the case of families of elliptic curves. The details will be discussed in a future paper.
\subsection{Contents}

In Section \ref{sec:Chen}, we give a quick and dirty exposition of Chen's reduced bar complex and his $\pi_1$-de Rham theorem for a complex manifold $M$. In the case $M = \mathbb{P}^1(\mathbb{C})\setminus \{0,1,\infty\}$, multiple zeta values make a natural appearance in Chen's theory; we briefly recall this in Section \ref{sec:P1}. Finally, in Section \ref{sec:E}, we describe the case of a once-punctured elliptic curve.

\subsection{Notation and conventions} \label{ssec:notationandconventions}
We will compose paths in the "algebraic geometer's order", i.e. $\gamma_1\gamma_2$ means to first travel along $\gamma_2$, then along $\gamma_1$.

Given a complex manifold $M$, we will denote by $\cA^{\bullet}(M)$ the $\bC$-valued \textit{$C^{\infty}$-de Rham complex of $M$}, viewed as a differential graded $\bC$-algebra. Given a smooth algebraic variety $X$ and a normal crossings divisor $D\subset X$, both defined over a field $k$, we will denote by $\Omega^{\bullet}_X\langle D\rangle$ the \textit{algebraic de Rham complex of $X$ with log poles along $D$}, \cite[\S 3.1]{Deligne:HodgeII}, \cite[\S 3.1.6]{HuberMuellerStach}.

For differential one-forms $\omega_1,\ldots,\omega_n \in \cA^1(M)$ and a piecewise smooth path $\gamma: [0,1] \rightarrow M$, we define the \textit{iterated integral} by
\begin{equation}
\int_\gamma\omega_1\ldots\omega_n:=\int_{1\geq t_1\geq \ldots\geq t_n\geq 0}\gamma^{\ast}(\omega_1)(t_1)\ldots\gamma^{\ast}(\omega_n)(t_n),
\end{equation}
which has \textit{length} $n$. For $n=0$, we adopt the convention $\int_\gamma \equiv 1$. Note that our convention for the order of iterated integration is compatible with our convention for path composition, in the sense that
\begin{equation}
\int_{\gamma_1\gamma_2}\omega_1\ldots\omega_n=\sum_{i=0}^n\int_{\gamma_1}\omega_1\ldots\omega_i\int_{\gamma_2}\omega_{i+1}\ldots\omega_n,
\end{equation}
for all piecewise smooth paths $\gamma_1,\gamma_2: [0,1]\rightarrow M$ with $\gamma_2(1)=\gamma_1(0)$. Given a point $p\in M$, we will often view iterated integrals as functions
\begin{equation}
\int\omega_1\ldots\omega_n: P(M) \rightarrow \bC, \quad \gamma \mapsto \int_\gamma\omega_1\ldots\omega_n,
\end{equation}
where $P(M)$ is the set of all piecewise smooth paths $\gamma:[0,1] \rightarrow M$. A $\bC$-linear combination of iterated integrals is called \textit{homotopy invariant} if for all points $p,q \in M$, its value at $\gamma \in P(M)$ with $\gamma(0)=p$, $\gamma(1)=q$, only depends on the homotopy class of $\gamma$ relative to $p$, $q$.

\section{Chen's reduced bar de Rham complex and the \texorpdfstring{$\pi_1$}{pi1}-de Rham theorem} \label{sec:Chen}

Throughout this section, we will denote by $M$ a connected complex manifold $M$ such that $H_1(M;\bQ)$ is finite-dimensional, for example $M=X(\bC)$ for a smooth algebraic variety $X$ over $\bC$.

We begin by giving some details for the reduced bar de Rham complex of $M$, \cite{Chen:ReducedBarComplex}. See also \cite[\S 7]{Hain:IteratedIntegrals} and \cite[\S 3.4.1]{BurgosFresan}. Let $A^{\bullet} \subset \cA^{\bullet}(M)$ be a differential graded $\bC$-subalgebra which is quasi-isomorphic to $\cA^{\bullet}(M)$. We will refer to such an $A^{\bullet}$ as a \textit{model} of $\cA^{\bullet}(M)$. For simplicity, we will assume furthermore that $A^{\bullet}$ is connected, that is $A^0\cong \bC$; see \cite{Chen:ReducedBarComplex} for the general case. In the connected case, the reduced bar complex $B^{\bullet}(A^{\bullet})$ of $A^{\bullet}$ is simply given by
$
B^{\bullet}(A^{\bullet}):=\bigoplus_{n\geq 0}(A^{>0}[1])^{\otimes n},
$
with differential
\begin{equation} \label{eqn:differential}
\begin{aligned}
d_B: B^{\bullet}(A^{\bullet})& \rightarrow B^{\bullet}(A^{\bullet})\\
[a_1|\ldots|a_n]& \mapsto \sum_{i=1}^n(-1)^i[Ja_1|\ldots|Ja_{i-1}|da_i|a_{i+1}|\ldots|a_n]\\
&\hphantom{+++++}+\sum_{i=1}^{n-1}(-1)^{i+1}[Ja_1|\ldots|Ja_{i-1}|a_i\wedge a_{i+1}|a_{i+2}|\ldots|a_n], 
\end{aligned}
\end{equation}
where $J: A^{\bullet}\rightarrow A^{\bullet}$ is given by $J(a)=(-1)^{\deg a}a$, and we write $[a_1|\ldots|a_n]$ instead of $a_1\otimes\ldots\otimes a_n$, as is customary. The reduced bar complex is concentrated in non-negative degrees and is equipped with an increasing filtration $B^{\bullet}_{\ast}(A^{\bullet})$ by length, where $B^{\bullet}_\ell(A^{\bullet}) \subset B^{\bullet}(A^{\bullet})$ denotes the subcomplex spanned by elements $[a_1|\ldots|a_n]$, with $n\leq \ell$. Note that $B^0(A^{\bullet})=\bigoplus_{n\geq 0}(A^1)^{\otimes n}$, so that 
\begin{equation}
H^0(B^{\bullet}(A^\bullet))=\Bigg\{ \xi \in \bigoplus_{n\geq 0}(A^1)^{\otimes n} \, \Big\vert \, d_B(\xi)=0 \Bigg\}.
\end{equation}
On the other hand, let $\bZ\pi_1(M,p)$ be the group ring of $\pi_1(M,p)$, and denote by $I=\ker(\bZ\pi_1(M,p) \stackrel{\gamma\mapsto 1}\longrightarrow \bZ)$ the augmentation ideal. There is a $\bC$-linear map
\begin{equation}
\begin{aligned}
\Int^\ell_p: H^0(B^{\bullet}_\ell(A^{\bullet})) &\rightarrow \Hom(\bZ\pi_1(M,p),\bC)\\
\sum_{n=0}^\ell\sum_{i_1,\ldots,i_n}[\omega_{i_1}|\ldots|\omega_{i_n}] &\mapsto \sum_{n=0}^\ell\sum_{i_1,\ldots,i_n}\int\omega_{i_1}\ldots\omega_{i_n}.
\end{aligned}
\end{equation}
which factors through $\Hom(\bZ\pi_1(M,p)/I^{\ell+1},\bC)$. For varying $\ell\geq 0$, the $\bC$-algebras $\Hom(\bZ\pi_1(M,p)/I^{\ell+1},\bC)$ form a filtered system whose transition maps 
\begin{equation}
\Hom(\bZ\pi_1(M,p)/I^{\ell+1},\bC) \rightarrow \Hom(\bZ\pi_1(M,p)/I^{\ell'+1},\bC), \quad \mbox{for }\ell'\geq \ell,
\end{equation}
are dual to the canonical projections $\bZ\pi_1(M,p)/I^{\ell'+1}\rightarrow \bZ\pi_1(M,p)/I^{\ell+1}$. The following result combines both Theorems \ref{thm:Chenintro} and \ref{thm:ChenintroII}
\begin{thm}[Chen] \label{thm:Chen}
For each $\ell\geq 0$, the induced map
\begin{equation}
\Int^\ell_p: H^0(B^{\bullet}_\ell(A^{\bullet})) \rightarrow \Hom(\bZ\pi_1(M,p)/I^{\ell+1},\bC)
\end{equation}
is an isomorphism.
\end{thm}
In fact, the union $H^0(B^{\bullet}(A^{\bullet}))=\bigcup_{\ell \geq 0}H^0(B^{\bullet}_\ell(A^{\bullet}))$ is a commutative filtered Hopf algebra with the shuffle product and deconcatenation coproduct, and the induced map $\Int_p$ is an isomorphism (of commutative filtered Hopf algebras). Also, note that the $\bQ$-algebra $\varinjlim_{\ell} \Hom(\bZ\pi_1(M,p)/I^{\ell+1},\bQ)$ is precisely the affine ring of functions on the pro-unipotent, or Mal'cev, completion of $\pi_1(M,p)$.

\section{The case \texorpdfstring{$\bP^1 \setminus \{0,1,\infty\}$}{P1}} \label{sec:P1}

We next look at Chen's $\pi_1$-de Rham theorem in the special case $M=\bP^1(\bC)\setminus \{0,1,\infty\}$, the complex projective line minus three points, which is the set of complex points of the algebraic curve $X:=\bP^1_\bQ\setminus \{0,1,\infty\}$.

\subsection{The logarithmic bar de Rham complex}
Let $A^{\bullet}_X=H^0(\overline{X},\Omega^{\bullet}_{\overline{X}}\langle D\rangle)$, which is given explicitly by
\begin{equation}
A^0_X=\bQ, \quad A^1_X=\bQ\omega_0\oplus\bQ\omega_1, \quad A^n_X=0, \quad n\geq 2,
\end{equation}
for $\omega_0=dz/z$, $\omega_1=dz/(z-1)$. Since all differentials and wedge products in $A^{\bullet}_X$ are trivial, the differential \eqref{eqn:differential} is trivial and the reduced bar complex of $A^{\bullet}_X$ is readily computed:
\begin{equation}
H^0(B^{\bullet}(A^{\bullet}_X)) \cong \bigoplus_{n\geq 0}(\bQ\omega_0\oplus\bQ\omega_1)^{\otimes n}.
\end{equation}
Equivalently, each iterated integral $\int\omega_{i_1}\ldots\omega_{i_n}$, where $i_j \in \{0,1\}$, is homotopy invariant.

Now by the remark just after \eqref{eqn:embedding}, $A^{\bullet}_X\otimes_\bQ\bC$ is a connected model for $\cA^{\bullet}(M)$, and for any base point $p \in X(\bQ)=\bP^1 \setminus \{0,1,\infty\}$, Theorem \ref{thm:Chen} implies that the $\bC$-linear map
\begin{equation} \label{eqn:Chenpi1}
\begin{aligned}
\Int_p: H^0(B^\bullet(A^{\bullet}_X)) \otimes_\bQ\bC &\rightarrow \varinjlim_\ell \Hom(\bZ\pi_1(M,p)/I^{\ell+1},\bC)\\
[\omega_{i_1}|\ldots|\omega_{i_n}] & \mapsto \int\omega_{i_1}\ldots\omega_{i_n}
\end{aligned}
\end{equation}
is an isomorphism. The image of $H^0(B^\bullet(A^{\bullet}_X))$ then defines a $\bQ$-structure on the right hand side which is given by the iterated integrals $\int_\gamma\omega_{i_1}\ldots\omega_{i_\ell}$, for $\gamma \in \pi_1(M,p)$. These complex numbers are examples of periods in the sense of Kontsevich--Zagier, \cite{KZ}.

\subsection{Chen's theorem for tangential base points}

In the definition of $\pi_1(M,p)$, the base point $p$ can be replaced by a nonzero tangent vector $\vec{v}_p \in T_p(\bP^1(\bC))$ where $p\in \{0,1,\infty\}$. A loop based at $\vec{v}_p$ is now required to satisfy $\gamma'(0)=\vec{v}_p$ and $\gamma'(1)=-\vec{v}_p$, and homotopies between two such paths should respect $\vec{v}_p$.

With these provisions, one can define as before the $\bC$-algebra $\Hom(\bZ\pi_1(M,b)/I^{\ell+1},\bC)$ where now $b$ is either point in $\bP^1(\bC)\setminus \{0,1,\infty\}$ or a non-zero tangent vector at one of the points $\{0,1,\infty\}$. The integration map $\Int_b$ relating the two objects is more involved if $b$ is a tangent vector, since the forms $\omega_0,\omega_1$ may have poles there. More specifically, one regularizes these iterated integrals using Deligne's formalism, \cite[\S 15.44]{Deligne:P1}, which essentially amounts to classical solution theory of ODEs with a regular singular point, and crucially uses that the forms $\omega_0,\omega_1$ have at most simple poles. In any case, we obtain a well-defined map
\begin{equation}
\Int_b: H^0(B^{\bullet}(A^{\bullet}_X)) \rightarrow \varinjlim_\ell \Hom(\bZ\pi_1(M,b)/I^{\ell+1},\bC)
\end{equation}
which induces an isomorphism of filtered $\bC$-vector spaces after tensoring the left hand side with $\bC$, \cite[Theorem 3.247]{BurgosFresan}.

We are particularly interested in the case $b=\pm \vec{1}_p$, for $p \in \{0,1,\infty\}$. In this case, the image of $H^0(B^{\bullet}(A^{\bullet}_X))$ under $\Int_b$ is uniquely determined by the values of the regularized iterated integrals $\int_\gamma\omega_{i_1}\ldots\omega_{i_\ell}$, and these are well known to evaluate to $\bQ[2\pi i]$-linear combinations of multiple zeta values
\begin{equation}
\zeta(k_1,\ldots,k_d):=\sum_{n_1>\ldots>n_d>0}\frac{1}{n_1^{k_1}\ldots n_d^{k_d}}, \quad k_1,\ldots,k_d\geq 1,\, k_1\geq 2.
\end{equation}
\begin{rmk}
With considerable work, one can show that for $b\in \{\pm\vec{1}_p \, \vert \, p \in \{0,1,\infty\}\}$ and every $\ell\geq 0$, the triple
\begin{equation}
(H^0(B^\bullet_\ell(A^{\bullet}_X)),\Hom(\bZ\pi_1(M,b)/I^{\ell+1},\bQ),\Int_b)
\end{equation}
is the de Rham--Betti realization of a mixed Tate motive over $\bZ$, \cite{DeligneGoncharov}. This interpretation leads to the definition of motivic multiple zeta values, \cite{Brown:Annals,Brown:SingleValuedPeriods}.
\end{rmk}
\section{The case of a once-punctured elliptic curve} \label{sec:E}

Let $E$ be an elliptic curve defined over a subfield $k\subset \mathbb{C}$, with origin $O\in E(k)$, and let $\tilde{E} = E\setminus\{O\}$ be the once-punctured elliptic curve. Let us also fix a Weierstrass model $\tilde{E} \cong \mathop{\rm Spec} k[x,y]/(y^2-4x^3+ax+b)$, for some $a,b \in k$ ($a^3-27b^2\neq 0$), and set $\omega = \frac{dx}{y} \in H^0(E,\Omega^1_E)$. Finally, we identify the complex manifold $\tilde{E}(\mathbb{C})$ with $U \coloneqq (\mathbb{C}/\Lambda)\setminus\{0\}$, where $\Lambda = \{\int_{\gamma}\omega \in \mathbb{C} \mid \gamma \in H_1(E(\mathbb{C}), \mathbb{Z})\}$.

Note that the differential graded $\mathbb{C}$-algebras $H^0(E, \Omega^{\bullet}_E\langle O \rangle)\otimes_k \mathbb{C}$ and $\mathcal{A}^{\bullet}(U)$  cannot be quasi-isomorphic:  on one hand, the Riemann-Roch formula implies that $H^0(E,\Omega_E^1\langle O \rangle) = H^0(E,\Omega^1_E) = k\omega$; on the other hand, $\dim_{\mathbb{C}} H^1(\mathcal{A}^{\bullet}(U)) = 2$. In order to get around this problem, we replace $E$ by its universal vectorial extension $\pi: E^{\dagger} \to E$, and likewise the complex $H^0(E,\Omega^{\bullet}_E\langle O \rangle)$ by  $H^0(E^{\dagger},\Omega^{\bullet}_{E^{\dagger}}\langle D \rangle)$, where $D \coloneqq \pi^{-1}(O)$.

\subsection{The universal vectorial extension of an elliptic curve}

See \cite[Appendix C]{Katz} for details. Let $E$ be an elliptic curve over $k$ with origin $O \in E(k)$ as above. We will identify $E$ in the standard way with $\Pic^0_{E/k}$, the moduli scheme of line bundles of degree zero on $E$. The universal vectorial extension $E^{\dagger}$ of $E$ is then defined to be the moduli scheme $\Pic^{\dagger}_{E/k}$ of pairs $(\cL,\nabla)$ where $\cL$ is a degree zero line bundle on $E$ and $\nabla$ is a $k$-linear integrable connection on $\cL$. This is a commutative algebraic $k$-group scheme which fits into a short exact sequence
\begin{equation}
0\longrightarrow \underline{\Omega}^1_E\longrightarrow E^{\dagger}\stackrel{\pi}\longrightarrow E\longrightarrow 0,
\end{equation}
where $\pi$ is the canonical projection induced by ``forgetting the connection'' and $\underline{\Omega}^1_{E}$ is the vector group defined by $\underline{\Omega}^1_E(R):=H^0(E_R,\Omega^1_{E_R})$, where $E_R:=E\otimes_k R$, for every $k$-algebra $R$. Using the generator  $\omega=\frac{dx}{y}$ of $H^0(E,\Omega^1_E)$, we may identify $\underline{\Omega}^1_{E} \cong \bG_a$ with coordinate $t$. Over the punctured curve $\widetilde{E}$, the projection $\pi$ splits canonically, \cite[\S C.2]{Katz}, giving rise to an isomorphism of schemes $\widetilde{E}^{\dagger}\cong \underline{\Omega}^1_{E} \times \widetilde{E}$, where $\widetilde{E}^{\dagger}:=E^{\dagger} \setminus D$, and $D:=\pi^{-1}(O)\cong \underline{\Omega}^1_{E}$.

In order to understand the structure of $E^{\dagger}(\bC)$ as a complex manifold, consider the subgroup $L:=\{(\lambda,-\eta(\lambda)) \, \vert \, \lambda \in \Lambda\} \subset \bC^2$, where $\Lambda$ is the lattice uniformizing $E(\bC)$, and $\eta(\lambda):=\int_\lambda\eta$, where $\eta=x\frac{dx}{y}$ is a differential of the second kind on $E$. In terms of the Weierstrass zeta function
\begin{equation}
\zeta(z)=\frac 1z+\sum_{\lambda \in \Lambda\setminus \{0\}}\left(\frac{1}{z-\lambda}+\frac{1}{\lambda}+\frac{z}{\lambda^2}\right)
\end{equation}
of the lattice $\Lambda$, we have $\eta(\lambda)=\zeta(z)-\zeta(z+\lambda)$, for any $z \in \bC\setminus \Lambda$. We then have a canonical isomorphism
\begin{equation} \label{eqn:uniformization}
\bC^2/L\cong E^{\dagger}(\bC),
\end{equation}
of complex manifolds, under which the divisor $D(\bC) \subset E^{\dagger}(\bC)$ corresponds to $\{0\}\times \bC \subset \bC^2/L$. On the open subset $\widetilde{E}^{\dagger}(\bC)$ with coordinate $(x,y,t)$, and writing $(z,s)$ for the canonical coordinate on $\bC^2$, the isomorphism \eqref{eqn:uniformization} is given explicitly by
\begin{equation}
\begin{aligned}
U^{\dagger}&\rightarrow \widetilde{E}^{\dagger}(\bC)\\
(z,s) &\mapsto (\wp(z),\wp'(z),\zeta(z)-s),
\end{aligned}
\end{equation}
where $U^{\dagger} \subset \bC^2/L$ is the complement of the divisor defined by $z=0$, and $\wp(z)=-\zeta'(z)$ is the Weierstrass $\wp$-function.
\begin{rmk}
The previous discussion implies in particular that we have an isomorphism $U^{\dagger}\cong U\times \bC$ of complex manifolds. Therefore, the canonical embedding
\begin{equation}
\cA^{\bullet}(U) \hookrightarrow \cA^{\bullet}(U^{\dagger}),
\end{equation}
of $C^{\infty}$-de Rham complexes, induced by pullback along $\pi$, is a quasi-isomorphism by the K\"unneth formula. In particular, the spaces of homotopy invariant iterated integrals on $U$ and on $U^{\dagger}$ are isomorphic.
\end{rmk}

\subsection{Differential forms on \texorpdfstring{$\widetilde{E}^{\dagger}$}{Edagger}}

While we have seen that the $C^{\infty}$-de Rham complexes of $U$ and $U^{\dagger}$ are quasi-isomorphic, it turns out that the complexes $H^0(E,\Omega^{\bullet}_E\langle O\rangle)$ and $H^0(E^{\dagger},\Omega^{\bullet}_E\langle D\rangle)$ of logarithmic forms are not. In particular, and unlike the former, the latter computes $H^*_{\rm dR}(U^{\dagger})$.

In order to compute $H^0(E^{\dagger},\Omega^1_{E^{\dagger}}\langle D\rangle)$, the basic fact is the following well known computation of the degree one Hodge cohomology groups.
\begin{prop} \label{prop:Hodgecohomology}
We have
\begin{equation}
H^1(E^{\dagger},\cO_{E^{\dagger}})=0, \quad H^0(E^{\dagger},\Omega^1_{E^{\dagger}})\cong H^1_{\rm dR}(E),
\end{equation}
where $H^1_{\rm dR}(E)$ denotes the first algebraic de Rham cohomology of $E$ over $k$.
\end{prop}
More precisely, there exists a $k$-basis $\omega$, $\nu$ for $H^0(E^{\dagger},\Omega^1_{E^{\dagger}})$ which on the open affine subset $\widetilde{E}^{\dagger}$ is given explicitly by $\omega=\frac{dx}{y}$, and $\nu=-dt-x\frac{dx}{y}$. Under the complex uniformization \eqref{eqn:uniformization}, the forms $\omega$, $\nu$ correspond to the one-forms $dz,ds \in \cA^1(U^{\dagger})$.

In order to write down a basis for $H^0(E^{\dagger},\Omega^1_{E^{\dagger}}\langle D\rangle)$, we consider the residue exact sequence of sheaves 
\begin{equation}
0 \longrightarrow \Omega^1_{E^{\dagger}} \longrightarrow \Omega^1_{E^{\dagger}}\langle D\rangle \stackrel{\rm Res}\longrightarrow \cO_D \longrightarrow 0.
\end{equation}
Combining this with Proposition \ref{prop:Hodgecohomology}, we get the following result (cf. \cite[Lemma 6.1]{EnriquezEtingof}).
\begin{prop}
We have an exact sequence of $k$-vector spaces
\begin{equation}
0 \longrightarrow H^1_{\rm dR}(E) \longrightarrow H^0(E^{\dagger},\Omega^1_{E^{\dagger}}\langle D\rangle) \stackrel{\rm Res}\longrightarrow k[t] \longrightarrow 0,
\end{equation}
where $k[t]\cong H^0(D,\cO_D)$.
\end{prop}
It is therefore enough to construct a splitting of the residue map. For this, we work on the complex uniformization $U^{\dagger}$ and show algebraicity only a posteriori. Define a family $\omega^{(n)}=f^{(n)}(z,s)dz$ of meromorphic one-forms on $\bC \times \bC$ via the generating series\footnote{We apologize for the aesthetically questionable act of putting the symbols $\omega$ and $w$ right next to each other.}
\begin{equation} \label{eqn:algebraic}
e^{-sw}\frac{\sigma(z+w)}{\sigma(z)\sigma(w)}dz=\sum_{n\geq 0}\omega^{(n)}w^{n-1},
\end{equation}
where $\sigma(z)=z\prod_{\lambda\in \Lambda\setminus \{0\}}(1-z/\lambda)e^{z/\lambda+z^2/2\lambda^2}$ is the Weierstrass sigma function. Clearly, we have $\omega^{(0)}=dz$.
\begin{prop} \label{prop:algebraic}
The one-forms $\omega^{(n)}$, for $n\geq 1$, descend to meromorphic one-forms on $\bC^2/L$ with simple poles along $z=0$ and residue $\frac{t^{n-1}}{(n-1)!}$. Under the uniformization \eqref{eqn:uniformization}, they correspond to global logarithmic one-forms on $\widetilde{E}^{\dagger}(\bC)$, which are defined over $k$. Finally, the forms $\omega^{(n)}$ satisfy the differential equation
\begin{equation} \label{eqn:differentialequation}
d\nu=d\omega^{(0)}=0, \quad d\omega^{(n)}=-\nu\wedge\omega^{(n-1)}, \quad n\geq 1.
\end{equation}
\end{prop}
\begin{rmk}
The forms $\omega^{(n)}$ are an algebraic variant of real analytic one-forms introduced by Brown--Levin, \cite[\S 3]{BrownLevin}, and have already appeared with slightly different conventions in \cite[Lemma 6.3]{EnriquezEtingof}.
\end{rmk}
By Proposition \ref{prop:algebraic}, the forms $\nu$, $\omega^{(n)}$, for $n\geq 0$, form a $k$-basis for $H^0(E^{\dagger},\Omega^1_{E^{\dagger}}\langle D\rangle)$ which splits the residue map. The space of two-forms $H^0(E^{\dagger},\Omega^2_{E^{\dagger}}\langle D\rangle)$ is determined from the following proposition.
\begin{prop}[{cf. \cite[Lemma 8]{BrownLevin}}] \label{prop:twoforms}
The family of two-forms $\{\nu\wedge\omega^{(n)}\}_{n\geq 0}$ is a $k$-basis of $H^0(E^{\dagger},\Omega^2_{E^{\dagger}}\langle D\rangle)$.
\end{prop}
\subsection{The logarithmic bar de Rham complex of \texorpdfstring{$E^{\dagger}$}{Edagger}}
Let $A^{\bullet}_{E^\dagger}:=H^0(E^{\dagger},\Omega^{\bullet}_{E^{\dagger}}\langle D\rangle)$. We have
\begin{equation}
A^0_{E^{\dagger}}=k, \quad A^1_{E^{\dagger}}=k\nu\oplus\bigoplus_{n\geq 0}k\omega^{(n)}, \quad A^2_{E^{\dagger}}=\bigoplus_{n\geq 0}k(\nu\wedge\omega^{(n)}),
\end{equation}
as well as $A^n_{E^{\dagger}}=0$, for $n\geq 3$. The following result is the analogue of \cite[Theorem 19]{BrownLevin}.
\begin{prop} \label{prop:quasiisomorphism}
The natural inclusion
$
A^{\bullet}_{E^{\dagger}} \otimes_k \bC \hookrightarrow \cA^{\bullet}(U^{\dagger}),
$
is a quasi-isomorphism.
\end{prop}
It remains to compute the reduced bar construction $H^0(B^{\bullet}(A^{\bullet}_{E^{\dagger}}))$. This is more involved than in the case of $\bP^1 \setminus \{0,1,\infty\}$ since $\widetilde{E}^{\dagger}$ is two-dimensional, so that integrability is now a non-trivial constraint. However, everything goes in almost the exact same way as in \cite{BrownLevin}, so we will simply state the results.

Let $\sx_0,\sx_1$ be formal variables, $\bL_k(\sx_0,\sx_1)$ be the free Lie $k$-algebra on $\{\sx_0,\sx_1\}$ and $\bL_k(\sx_0,\sx_1)^{\wedge}$ its completion for the lower central series. Consider the formal differential one-form\footnote{The notation is chosen because $\omega_{\rm KZB}$ is essentially an algebraic version of the elliptic Knizhnik--Zamolodchikov--Bernard connection form, \cite{CEE,Hain:KZB,LevinRacinet}.}
\begin{equation}
\omega_{\rm KZB}:=\nu\otimes \sx_0+\sum_{n\geq 0}\omega^{(n)}\otimes\ad_{\sx_0}^n(\sx_1) \in A^1_{E^{\dagger}}\widehat{\otimes}\bL_k(\sx_0,\sx_1)^{\wedge}, \quad \ad^n_{\sx_0}(\sx_1):=[\sx_0,\ad_{\sx_0}^{n-1}(\sx_1)]
\end{equation}
where all tensor products are over $k$. It follows from \eqref{eqn:differentialequation} that $\omega_{\rm KZB}$ is integrable:
\begin{equation}
d\omega_{\rm KZB}+\omega_{\rm KZB}\wedge\omega_{\rm KZB}=0.
\end{equation}
Now consider the formal series
\begin{equation}
I(\omega_{\rm KZB})=\sum_{\ell \geq 0}\underbrace{[\omega_{\rm KZB}|\ldots|\omega_{\rm KZB}]}_{\ell} \in \bigoplus_{\ell \geq 0}(A^1_{E^{\dagger}})^{\otimes \ell} \widehat{\otimes} \, \widehat{\cU}(\bL_k(\sx_0,\sx_1)^{\wedge}),
\end{equation}
where $\widehat{\cU}(\bL_k(\sx_0,\sx_1)^{\wedge}) \cong k\langle\!\langle \sx_0,\sx_1\rangle\!\rangle$ denotes the completed universal enveloping algebra of $\bL_k(\sx_0,\sx_1)^{\wedge}$. By duality, each word $w\in \langle \sx_0,\sx_1\rangle$ in the letters $\sx_0,\sx_1$ defines a map
\begin{equation}
\begin{aligned}
c_w: \bigoplus_{\ell \geq 0}(A^1_{E^{\dagger}})^{\otimes \ell} \widehat{\otimes} \, \widehat{\cU}(\bL_k(\sx_0,\sx_1)^{\wedge}) &\rightarrow \bigoplus_{\ell \geq 0}(A^1_{E^{\dagger}})^{\otimes \ell},\\
v&\mapsto \delta_{v,w},
\end{aligned}
\end{equation}
which simply singles out the coefficient of $w$. The following proposition is the analogue of \cite[Proposition 23]{BrownLevin}.
\begin{prop}
We have
\begin{equation}
H^0(B^{\bullet}(A^{\bullet}_{E^{\dagger}})) \cong \spn_k\{c_w(I(\omega_{\rm KZB})) \, \vert \, w \in \langle \sx_0,\sx_1\rangle \}.
\end{equation}
\end{prop}
Applying Chen's theorem, we obtain an isomorphism of $\bC$-vector spaces
\begin{equation}
\Int_p: H^0(B^{\bullet}(A^{\bullet}_{E^{\dagger}})) \otimes_k\bC \rightarrow \varinjlim_\ell \Hom(\bZ\pi_1(U^{\dagger},p)/I^{\ell+1},\bC).
\end{equation}
Moreover, since $H^0(B^{\bullet}(A^{\bullet}_{E^{\dagger}}))$ is defined over $k$, its image gives a $k$-structure on the right hand side, which is given by the values of $c_w(I(\omega_{\rm KZB}))$ integrated along loops on $U^{\dagger}$ based at $p$. In the case where $p \in \widetilde{E}^{\dagger}(k)$ and $k \subset \bC$ is finite over $\bQ$, these iterated integrals are periods in the sense of Kontsevich--Zagier. 
\begin{exmp}
In the simplest case $\ell=1$, the $k$-structure is given by the (quasi-) periods of elliptic curves. To see this, let $E$ be an elliptic curve over $k$ with complex uniformization $E(\bC) \cong \bC/\Lambda$, and consider the length one elements $[\omega^{(0)}]$, $[\nu] \in H^0(B_1^{\bullet}(A^{\bullet}_{E^{\dagger}}))$. Applying Chen's isomorphism, these correspond to the line integrals
\begin{equation}
\int\omega^{(0)}, \, \int\nu \in \Hom(\bZ\pi_1(U^{\dagger},p)/I^2,\bC) \cong \bC\oplus \Hom(\Lambda,\bC),
\end{equation}
which are independent of the choice of base point $p$. Recalling that the one-forms $\omega^{(0)}$, $\nu$ pullback to the one-forms $dz$, respectively $ds$, on $T^{\dagger}\cong \bC^2/L$, the functionals $\int\omega^{(0)}, \, \int\nu: H_1(E;\bZ) \rightarrow \bC$ are then given by
\begin{equation}
\int_{\lambda}\omega^{(0)}=\lambda, \quad \int_{\lambda}\nu=-\eta(\lambda).
\end{equation}
\end{exmp}

\bigskip

\noindent
{\bf Acknowledgements:}
Many thanks to Richard Hain for his valuable comments and corrections on an earlier version of the manuscript. The second author wishes to thank both Hidekazu Furusho and Masanobu Kaneko, as well as the Research Institute for Mathematical Sciences (RIMS) and Kyushu University, for hospitality. This project has received funding from the European Research Council (ERC) under the European Union’s Horizon 2020 research and innovation programme (grant agreement No. 724638)

\end{document}